
\input amstex
\documentstyle{amsppt}
\NoBlackBoxes

\magnification=\magstep1

\topmatter

\font\midbf=cmbx10 scaled1200

\def\QED{\hfill{\bf Q.E.D.}}
\def\lra{\longrightarrow}
\def\wt{\widetilde}

\title  $ $ \\ \\
    THE~MINIMAL~NUMBER~OF~SINGULAR~FIBERS
    \\    OF~A~SEMISTABLE~CURVE~OVER~${\Bbb P^1}$
\endtitle

\author{ SHENG-LI TAN }\endauthor

\leftheadtext\nofrills{SHENG-LI TAN}
\rightheadtext\nofrills{MINIMAL NUMBER OF SINGULAR FIBERS}

\address
 Department of Mathematics,
 East China Normal University,
 Shanghai 200062,
 P.~R.~of China\endaddress

\curraddr
 Max-Planck-Institut f{\"u}r Mathematik,
 Gottfried-Claren-Str. 26,
 53225 Bonn, Germany \endcurraddr

\email tan\@mpim-bonn.mpg.de \endemail

\thanks
The author would like to thank the hospitality and financial support of
Max-Planck-Institut f{\"u}r Mathematik in Bonn during this research.
The author thanks also Universit\'e de Nice for inviting me to give
a talk on this work in the Seminar of Projective Geometry.
This research is partially supported by the National Natural Science
Foundation of China and by the Science Foundation of the University
Doctoral Program of CNEC.
\endthanks

\endtopmatter

\document
\baselineskip15pt

\vskip0.5cm
\head{\bf 1. Introduction}
\endhead

   The purpose of this paper is to try to answer the following

\proclaim{\indent Szpiro's Question } {\rm ([Sz],  [B1]) Let $f:
S\longrightarrow
\Bbb P^1_{\Bbb C}$ be a family of
semistable curves of genus $g$, which is not trivial.
Then, what is the minimal number of the singular fibers of $f$?
}\endproclaim

   Beauville gives a lower bound for the number of
singular fibers.

\proclaim{\indent Beauville's Theorem} {\rm [B1]} With the notations as above,
if $g\geq1$,
then

 {\rm 1)} $f$ admits at least 4 singular fibers.

 {\rm 2)} If $f$ admits 4 singular fibers, then $S$ is algebraically
simply connected with $p_g(S)=0$, and the irreducible components of the 4
singular fibers are rational curves {\rm(}may be singular{\rm)}, which
generate a hyperplane of the $\Bbb Q$-vector space
{\rm Pic}$(S)\otimes \Bbb Q$.
\endproclaim

   Furthermore, Beauville [B1] gives an example of semistable elliptic
fibration over
$\Bbb P^1$ with 4 singular fibers, and one example of genus 3 with 5 singular
fibers, and  he gives also a series of such examples
with 6 singular fibers for all $g>1$.
In fact, Beauville conjectured
that for $g \geq 2$, there is no such fibrations with 4 singular
fibers.  In [B2], Beauville classified all semistable elliptic fibrations
over $\Bbb P^1$ with 4 singular fibers.

    Szpiro [Sz] considered that problem over a field with characteristic $p>
0$,
and he proved that the minimal number of singular fibers is at least 3, and if
 the surface is of general type, then the number is at least 4.

   The main result of this paper is

\proclaim{\indent Theorem 1 (Beauville's conjecture)} If $f: S\longrightarrow
\Bbb P^1_{\Bbb C}$
 is a non-trivial
semi-stable fibration of genus $g \geq 2$, then
 $f$ admits at least 5 singular fibers.
\endproclaim

Theorem~1 is an immediate consequence of Beauville's theorem
and the following ``strict canonical class inequality'' (cf. Sect. 2.1).

\proclaim{\indent Theorem~2} Let $f:S\longrightarrow  C$ be a locally
non-trivial semistable fibration of genus $g\geq2$ with $s$ singular fibers.
Then we have
$$
     \deg f_*\omega_{S/C}< {g\over 2}(2g(C)-2+s).
$$
\endproclaim

  The validity of these two theorems is heavily dependent on Miyaoka-Yau
inequality.
In Sect.~4, we shall give an example $f:S\longrightarrow  \Bbb P^1$ of genus 2
with 5
singular fibers. Note that Beauville has also given such an example for $g= 3$.

  In the first version of this paper, I proved Beauville's conjecture
for $g\geq 6$ by using Miyaoka's inequality. In the second version (MPI
preprint 94-45), the two theorems of this paper are proved by considering
Kodaira-Parshin's construction.
After the second version submitted, I simplified the proof and used the method
to obtain a linear (in $g$) and effective height inequality for algebraic
points
on curves over functional fields [Ta]. Later, the editor informed me that,
based
on the first version of this paper,
Nguyen Khac Viet had also
independently given the simple proof of Beauville's conjecture.

   I would like to thank Prof. A. Beauville, Prof. F. Hirzebruch, and Prof. G.
Xiao
for their helps and encouragements.  Prof. Xiao kindly informed me that he had
independently
obtained most of the steps in orignal version except for the proof of the main
theorems.
Special thanks are due to the referee for the suggestion for correction of the
original
version.

\vskip1cm
\head{\bf 2. Preliminaries}
\endhead

{\bf 2.1. Invariants of fibrations.}
\
   Let $f:S \longrightarrow  C$ be a relatively minimal fibration of genus $g
\geq2$, i.e., $S$ contains no $(-1)$-curves in a fiber of $f$. Let

   $$ \eqalign{&\chi_f=\deg f_*\omega_{S/C}=\chi(\Cal O_S)-(g-1)(g(C)-1),\cr
               &K_{S/C}^2=K_{S}^2-8(g-1)(g(C)-1),\cr
               &e_f=\chi_{\text{top}}(S)-4(g-1)(g(C)-1).\cr}
   $$
\noindent
   They are the basic invariants of $f$. If $K_{S/C}$ is the relative
canonical divisor of $f$, then $\chi_f=\deg f_*\omega_{S/C}$,
where $\omega_{S/C}=\Cal O(K_{S/C})$.
If $f$ is not locally trivial, by the well-known Arakelov-Parshin
Theorem ([Ar], [Pa]), we know $\chi_f>0$ and $K_{S/C}^2>0$. Then we can define
  the slope of $f$ as $\lambda_f={K_{S/C}^2/ \chi_f}$,
 which is an important invariant of $f$.
In [Xi], Xiao shows that if $f$ is a locally non-trivial fibration of genus
$g\geq 2$,
then we have
      $$\lambda_f\geq 4-{4\over g}.
                                                                   \eqno(1) $$
Furthermore, if the slope of $f$ is $4-{4/g}$,
then $e_f>0$, hence $f$ admits at least one singular fiber. Cornalba and Harris
[CH] have also obtained (1) for semistable fibrations.

{\bf 2.2. Miyaoka's inequality and Vojta's inequality.}
\
We refer to [Hi] for the details of the following Miyaoka's inequality.
\proclaim{\indent Lemma 2.3} {\rm [Mi]} If $S$ is a smooth surface such that
the canonical divisor $K_S$ is nef (numerically effective),
  and
$E_1$, $\cdots$, $E_n$ are disjoint $ADE$ curves on $S$, then we have

     $$
       \sum_{i=1}^nm(E_i)\leq 3c_2(S)-c_1^2(S),
     $$
where $m(E)$ is defined as follows,
     $$
       \eqalign{ m(A_r)&=3(r+1)-{3\over r+1} ,\cr
                 m(D_r)&=3(r+1)-{3\over 4(r-2)}, \hskip0.2cm \text{ for } r\geq
4,                     \cr
                 m(E_6)&=21-{1\over8},                     \cr
                 m(E_7)&=24-{1\over16},                     \cr
                 m(E_8)&=27-{1\over40}.                     \cr}
     $$
\endproclaim

    Finally, we should mention  Vojta's interisting inequality
 for semistable fibrations $f$ with $s$ singular
fibers, i.e., the ``canonical class inequality'' [Vo]:

    $$K_{S/C}^2\leq (2g-2)(2g(C)-2+s).
\eqno(2)$$
Our proof of Theorem 2 is based on these inequalities and Beauville's theorem.

\vskip1cm
\head{\bf 3. The proof of Theorem~2}
\endhead

First of all, we give some notations. Let $f:S\lra C$ be a relatively minimal
semistable fibration with $s$ singular fibers.
We denote by $f^\#: S^\#\lra C$ the corresponding stable model, and by $q$ a
singular point of $S^\#$.
Then $q$ is a rational double point of type
 $A_n$. Let $\mu_q=n$ be the Milnor
number of $(S^\#,q)$, i.e., the number of $(-2)$-curves  in the exceptional set
$E_q$ of the minimal resolution of $q$. We also denote by $q$
a singular point of a fiber on the smooth part of $S^\#$, in this case $\mu_q
=0$.

\proclaim{\indent Lemma 3.1} If $s>0$, then
    $$K_{S/C}^2< (2g-2)(2b-2+s). \eqno(3)$$
\endproclaim
\demo{\indent Proof}
  First we consider the base changes $\pi: \widetilde C \longrightarrow  C$ of
degree $de$
ramified exactly over the $s$ critical points,
where $d$ and $e>1$ are natural numbers. We assume that the fiber of $\pi$ over
each critical
point consists of $d$ points with ramification index $e$.  By Kodaira-Parshin
construction such a base change exists for all $e$ if $b>0$ and for odd $e$ if
$b=0$. (cf. [Vo]). Let $\widetilde f: \widetilde S\longrightarrow  \widetilde
C$ be the
pullback fibration of $f$ under $\pi$, then it is easy to know that
the inequality (3) holds if and only if
it holds for $\widetilde f$. By Beauville's theorem, if $b=0$, then $s\geq 4$,
so
$g(\widetilde C)>0$.
Hence we can assume that $b> 0$. Now we know that
    $$K_S\sim K_{S/C}+(2b-2)F$$
is nef, where $F$ is a fiber of $f$. Since $f$ is semistable,  we know that
$e_f$ is the number of singular points of the fibers, which implies
$e_f=\sum_q(\mu_q+1)$. Then by using Miyaoka's inequality,
we can obtain easily that
    $$
     K_{S/C}^2\leq \sum_q{3\over \mu_q +1}+(2g-2)(2b-2).       \eqno(4)
    $$
(cf. Vojta's proof. [Vo]).

Consider the pullback fibration $\wt f$ mentioned above, and denote by $\wt
\cdot$ the
corresponding objects of $\wt f$,  we have
    $$\eqalign{
\wt s=ds, \hskip0.3cm &K^2_{\wt S/\wt C}=deK^2_{S/C},
\hskip0.3cm \mu_{\wt q}+1=e(\mu_q+1),\cr
  2\wt b&-2=de(2b-2)+d(e-1)s.}
    $$
Applying (4) to $\wt f$, we have
     $$deK^2_{S/C}\leq {d\over e}\sum_q{3\over
\mu_q+1}+(2g-2)((2b-2)de+d(e-1)s),
     $$
i.e.,
    $$K_{S/C}^2-(2g-2)(2b-2+s)\leq -{(2g-2)s\over e} +{1\over e^2}\sum_q{3\over
       \mu_q+1}, \eqno(5)
     $$
let $e$ be large enough we can see that the right hand side of (5) is negative.
This completes
the proof of the lemma.
\QED
\enddemo

\vskip1cm
{\it Proof of Theorem 2}. \
{}From (1) and (2), we have
$$\chi_f \leq {g\over 2}(2b-2+s).$$
If the above equality holds, then
the equalities in (1) and (2) hold.
On the other hand, a locally non-trivial fibration with $\lambda_f=4-4/g$
admits at least one singular fiber, hence the equality in (2) contradicts Lemma
3.1. This completes the proof of Theorem 2.
\QED

\vskip1cm
\head{\bf 4. An example of genus $g=2$ with $s=5$}
\endhead

In this section, we shall construct a semistable fibration $f:S\longrightarrow
\Bbb P^1$
of genus 2 with 5 singular fibers.

Let $\phi$ and $\psi: \Bbb P^1\longrightarrow  \Bbb P^1$ be two morphisms with
$\deg \phi+\deg \psi=2g+2$. We assume that there exists a subset
$R=\{p_1,\cdots, p_5\}\subset \Bbb P^1$ satisfying

    i) the branched points of $\phi$ and $\psi$ are contained in $R$, and the
ramification points of them are of index 2.

    ii) $\phi^{-1}(p_i)\cap\psi^{-1}(p_i)$ consists of non-ramified points of
$\phi$ and $\psi$, and if $p\not\in R$, then $\phi^{-1}(p)\cap \psi^{-1}(p)$ is
empty.

In $\Bbb P^1\times \Bbb P^1$, we consider the divisors $\Gamma_{\phi}$ and
$\Gamma_{\psi}$, graphes of $\phi$ and $\psi$ respectively. Let
$B=\Gamma_{\phi}+
\Gamma_{\psi}$. Then $B$ is an even divisor of type $(2g+2,2)$.
Let  $\pi:\Sigma\longrightarrow  \Bbb P^1\times\Bbb P^1$ be
a double cover branched along $B$, and let $S$ be the canonical resolution
of the singularities of $\Sigma$. Then the second projection
$\Bbb P^1\times \Bbb P^1\longrightarrow  \Bbb P^1$ induces a semistable
fibration of
genus $g$ with 5 singular fibers.

Now we give an example of genus 2 with 5 singular fibers.
Let $a$ and $b$ be two nonzero complex numbers such that the discriminant of
the polynomial
    $$
p(x)=x^3+(2a-b^2)x^2+(a^2+2ab^2)x-a^2b^2,
    $$
is zero. Hence $p(x)$ has (at most) two zeros $x_1$ and $x_2$. Let $\phi$ and
$\psi: \Bbb P^1\longrightarrow  \Bbb P^1$ be two morphisms defined by
    $$
      \phi(t)=t^2+{a^2\over t^2}, \hskip0.5cm
      \psi(t)=-2a{t^2+b^2\over t^2-b^2}.
    $$
Note that if $\phi(t)=\psi(t)$, then $p(t^2)=0$. Take
$R=\{ \infty, 2a, -2a, x_1+a^2/x_1,x_2+a^2/x_2\}$. It is easy to check that
$\phi$ and $\psi$ satisfy i) and ii). This completes the construction.

\enddocument